\documentclass[11pt]{amsart}

\usepackage{amsfonts}
\usepackage{amssymb}
\usepackage{amscd}
\input xy
\xyoption{all}

\def\ZZ{{\mathbb Z}}

\def\QQ{{\mathbb Q}}
\def\PP{{\textbf P}}
\def\OO{{\mathcal O}}

\def\F{\mathcal{F}}
\def\E{\mathcal{E}}
\def\G{\mathcal{G}}

\def\I{\mathcal{I}}

\def\Pic0{{\rm Pic}^0(X)}

\theoremstyle{plain}

\newtheorem{theorem}{Theorem}[section]
\newtheorem{proposition}[theorem]{Proposition}
\newtheorem{corollary}[theorem]{Corollary}
\newtheorem{lemma}[theorem]{Lemma}

\theoremstyle{definition}
\newtheorem{definition}[theorem]{Definition}
\newtheorem{remark}[theorem]{Remark}
\newtheorem{example}[theorem]{Example}

\newtheorem{conjecture/question}[theorem]{Conjecture/Question}
\newtheorem{question}[theorem]{Question}
\newtheorem{remark/definition}[theorem]{Remark/Definition}

\theoremstyle{remark}

\begin{document}

\title{Regularity on abelian varieties III: further applications}

\author[G. Pareschi]{Giuseppe Pareschi}
\address{Dipartamento di Matematica, Universit\`a di Roma, Tor Vergata, V.le della
Ricerca Scientifica, I-00133 Roma, Italy}
\email{{\tt pareschi@mat.uniroma2.it}}

\author[M. Popa]{Mihnea Popa}
\address{Department of Mathematics, Harvard University,
One Oxford Street, Cambridge, MA 02138, USA }
\email{{\tt mpopa@math.harvard.edu}}

\thanks{The second author was partially supported by the NSF grant DMS-0200150.}

\date{\today}
\maketitle

\tableofcontents

\setlength{\parskip}{.1 in}

\markboth{G. PARESCHI and M. POPA}
{\bf REGULARITY ON ABELIAN VARIETIES III}

\section{Introduction}

Recently we have developed a regularity theory for coherent sheaves
on abelian varieties, called \emph{$M$-regularity} (cf. \cite{pp1}, \cite{pp2}).
It is a technique geared (at the moment) mainly towards solving geometric
problems related to linear series or equations for (subvarieties of)
abelian varieties. The main ingredients are the derived category
theoretic context of the Fourier-Mukai functor and the systematic use
of cohomological techniques. We refer the reader to the above
mentioned papers for full details. The more restricted purpose of the present sequel
is to describe a number of new applications of this theory in several
different directions in the study of abelian varieties and irregular varieties, some of which
we have announced in the previous papers.

We start in \S2 by giving an overview of the general context of
$M$-regularity. This essentially amounts to
recalling some basic definitions and results from \cite{pp1} which
will be used in the subsequent sections.

In \S3 we address a very familiar problem in the context of effective
results for linear series, namely that of bounding the Seshadri
constant measuring the local positivity of an ample line bundle.
We refer the reader to the upcoming book \cite{rob2} for a
comprehensive survey of the main results in this area.
A considerable body of work on this problem has developed in the context
of abelian varieties, where the Seshadri constants turn out to have
interesting connections with metric or arithmetic invariants
(cf. \cite{rob1}, \cite{nakamaye}, \cite{bauer1}, \cite{bauer2},
\cite{debarre} and
also \cite{rob2} for further references). Here we
explain how the Seshadri constant of a polarization $L$ on an abelian
variety is bounded below by an asymptotic version -- and in particular
by the usual -- $M$-regularity index of the line bundle $L$, as
defined in \cite{pp2}. This is the content of Theorem \ref{seshadri}.
The proof is a simple application of the $M$-regularity criterion of
\cite{pp1}, via the techniques of \cite{pp2} \S3. Combining this with
various bounds for Seshadri constants proved in \cite{rob1}, we obtain
bounds for $M$-regularity indices which are not apparent otherwise.
A problem of more interest -- at least historically -- and for which
we do not have an answer at this stage, is to produce uniform bounds on Seshadri constants
by bounding directly the $M$-regularity indices and applying the
result we prove here. We raise a few questions in this direction.

In \S4 we shift our attention towards a cohomological study of
Picard bundles. These are vector bundles on the Jacobian of a
curve $C$, whose projectivization is the symmetric product $C_n$,
for large enough $n$ -- in other words they parametrize all linear
series of degree $n$ on the curve $C$. Picard bundles have  been
the focus of intensive study, especially since they are closely
related to Brill-Noether theory; cf. \cite{rob2} 6.3.C and 7.2.C
for a general discussion and the corresponding literature. (It
seems in fact that the very definition of the Fourier functor was
given by Mukai in \cite{mukai1} in part with the aim of studying
Picard bundles.) We combine Fourier-Mukai techniques with the use
of the Eagon-Northcott resolution for special determinantal
varieties in order to compute their (strong) Theta regularity. In
down to earth terms this amounts to the following: it is known
that all Picard bundles are negative (i.e. have ample dual
bundle). However, we show that as soon as we twist them with the
smallest possible polarizations, namely the theta divisor and all
its translates, their higher cohomology vanishes. The same holds
for all their relatively small tensor powers (cf. Theorem
\ref{mreg_picard}). This vanishing theorem has numerous practical
applications, as does basically any nontrivial statement on Picard
bundles. In particular we recover in a more direct fashion the
main results of \cite{pp1} \S4 on the equations of the $W_d$'s  in
Jacobians, and on vanishing for pull-backs of pluritheta line
bundles to symmetric products. We end the section with a new
result on the equations of $Sing(\Theta)$ on non-hyperellptic
jacobians.

In \S5 we approach the problem of giving effective results for
linear series on irregular varieties of maximal Albanese dimension
via $M$-regularity for direct images of canonical (or adjoint)
bundles, extending work in \cite{pp1} \S5. The main result is a
theorem to the effect that, on a smooth projective minimal variety
$Y$ of general type whose Albanese map is generically finite and
whose Albanese image is not ruled by subtori, the pluricanonical
map given by $|\omega_Y^{\otimes 3}|$ is very ample outside the
locus of non-finite fibers. A result of this type was also proved
by Chen and Hacon \cite{chenhacon}, under a slightly more general
hypothesis, but with a less explicit conclusion. The relevance of
such statements, especially the fact that the hypothesis is not
too restrictive, comes from results of Green, Ein and Lazarsfeld
\cite{greenlaz}, \cite{einlaz} in the context of generic vanishing
theorems, explained in \S5. As an example \cite{einlaz}, Theorem
3, says that for any variety of maximal Albanese dimension,
$\chi(\omega_Y) = 0$ implies that the Albanese image of $Y$ is
ruled by tori. We note that at least part of this is purely
conceptually explained by the notion of $M$-regularity. We also
observe that results similar to the theorem stated above can be
obtained for higher order jets, and especially also for pluri-adjoint linear 
series on any variety of maximal Albanese dimension.

In \S6 we concentrate on the study of higher rank vector
bundles on abelian varieties. By work of Mukai and others (\cite{mukai3}, \cite{mukai4}, \cite{mukai1},
\cite{umemura} and \cite{orlov}) it has emerged that on abelian varieties the class of vector bundles
most closely resembling semistable vector bundles on curves and line
bundles on abelian varieties is that of \emph{semihomogeneous} vector bundles.
A vector bundle $E$ is semihomogeneous if every translation $t_x^* E$
by an element in $X$ is isomorphic to a twist $E\otimes \alpha$ by a
line bundle $\alpha\in {\rm Pic}^0 (X)$. It turns out that these
bundles are semistable, behave nicely under isogenies and
Fourier-Mukai transforms, and have a Mumford type theta-group theory
as in the case of line bundles. All of these suggest that there should
exist numerical criteria for their geometric properties like global or
normal generation. We show here that this is indeed the case, and the
measure is precisely the Theta regularity of the bundles in
question. More generally, we give a result on the surjectivity of the
multiplication map on global sections for two such vector bundles
(cf. Theorem \ref{mixed_multiplication}). Basic examples are the projective normality of
ample line bundles on any abelian variety, and the normal generation
of the Verlinde bundles on the Jacobian of a curve, coming from moduli
spaces of vector bundles on that curve. Note again that, although this
was part of our initial motivation, we do not have to appeal to the
theta-group theory of semihomogeneous bundles.

\noindent
\textbf{Acknowledgements.} We would like to thank Rob Lazarsfeld for
having introduced us to some of these topics and for interesting suggestions.
We also thank Christopher Hacon for discussions, and Olivier Debarre
for pointing out a serious mistake in \S4.

\section{$M$-regularity for coherent sheaves}

In this section we recall the main definitions and results from \cite{pp1}.
Let $X$ be an abelian variety of dimension $g$ over an algebraically
closed field. Given a coherent sheaf $\F$ on $X$, we denote
$S^i(\F):={\rm Supp}(R^i\hat{\mathcal S}(\F))$.
The sheaf $\F$ on $X$ is called \emph{$M$-regular} if
${\rm codim}(S^i(\F))>i$ for any $i=1,\ldots ,g$ (where, for $i=g$,
this means that $S^g(\F)$ is empty). If ${\rm Supp}(R^i\hat{\mathcal S}(\F))=\emptyset$,
this condition is trivially verified, and the sheaf is said to satisfy the Index Theorem
(I.T.) with index $0$. By the base change theorem, this is equivalent to saying that
$H^{i}(\F\otimes\alpha)=0$ for all $\alpha\in {\rm Pic}^{0}(X)$ and all $i>0$.
Finally, an extremely useful concept in the context of irregular varieties is the following:

\begin{definition}(\cite{pp1}, Definition 2.10)\label{continuous}
Let $Y$ be an irregular variety. We define a sheaf $\F$ on $Y$ to
be \emph{continuously globally generated} if for any non-empty open subset
$U\subset {\rm Pic}^{0}(Y)$ the sum of evaluation maps
$$\bigoplus_{\alpha\in U} H^0(\F\otimes \alpha)\otimes \alpha^\vee \longrightarrow \F$$
is surjective.
\end{definition}

\begin{theorem}({\bf $M$-regularity criterion},
\cite{pp1} Theorem 2.4 and Proposition 2.13.)\label{M_reg} Let
$\F$ be an $M$-regular sheaf on $X$, possibly supported on a
subvariety $Y$ of $X$. Then the following hold:
\newline
\noindent
(a) $F$ is continuously globally generated.
\newline
\noindent
(b) Let also $A$ be a line bundle on $Y$, continuously globally generated
as a sheaf on $X$.
Then $F\otimes A$ is globally generated.
\end{theorem}

\begin{theorem}(\cite{pp1} Theorem 2.5)\label{surjectivity}
Let $\F$ and $H$ be sheaves on $X$ such that $\F$ is $M$-regular and
$H$ is locally free satisfying
I.T. with index $0$. Then, for any non-empty Zariski open set $U\subset \hat X$, the map
$${\mathcal M}_U:\bigoplus_{\xi\in U}H^0(X,\F\otimes P_\xi)\otimes H^0(X,H\otimes P_\xi^\vee)
\buildrel{\oplus m_\xi}\over
\longrightarrow H^0(X,\F\otimes H)$$
is surjective, where $m_{\xi}$ denote the multiplication maps on global sections.
\end{theorem}

\begin{definition}
A coherent sheaf $\F$ on a polarized abelian variety $(X, \Theta)$ is
called $m$-$\Theta$-\emph{regular} if $\F((m-1)\Theta)$ is $M$-regular.
We will call it \emph{strongly} $m$-$\Theta$-\emph{regular} if $\F((m-1)\Theta)$
satisfies I.T. with index $0$.
\end{definition}

We recall the "abelian" Castelnuovo-Mumford Lemma, which is in fact a consequence
of the two results above.

\begin{theorem}(\cite{pp1} Theorem 6.3)\label{acm}
Let $\F$ be a $0$-$\Theta$-regular coherent sheaf on $X$. Then:
\newline
\noindent
(1) $\F$ is globally generated.
\newline
\noindent
(2) $\F$ is $m$-$\Theta$-regular for any $m\geq 1$.
\newline
\noindent
(3) The multiplication map
$$H^0(\F(\Theta))\otimes H^0(\OO(k\Theta))\longrightarrow H^0(\F((k+1)\Theta))$$
is surjective for any $k\geq 2$.
\end{theorem}

\section{$M$-regularity indices and Seshadri constants}

Here we express a natural relationship between Seshadri constants of
ample line bundles on abelian varieties and the $M$-regularity indices
of those line bundles as defined in \cite{pp2}. This result improves the
lower bound for Seshadri constants proved in \cite{nakamaye}, and
combined with the results of \cite{rob1} provides bounds for
controlling $M$-regularity. For a general overview of Seshadri
constants -- and in particular the statments used below -- one
can consult \cite{rob2} Ch.I \S5.

We start by recalling the basic definition from \cite{pp2} and by also
looking at a slight variation. We will denote by $X$ an abelian variety of
dimension $g$ over an algebraically closed field and by $L$ an ample
line bundle on $X$.

\begin{definition}
The $M$-\emph{regularity index} of $L$ is defined as
$$m(L):={\rm max}\{l~|~L\otimes
m_{x_1}^{k_1}\otimes \ldots \otimes m_{x_p}^{k_p}~{\rm is ~}M{\rm -regular~ for
~all~distinct~}$$
$$x_1,\ldots,x_p\in X {\rm~with~} \Sigma k_i=l\}.$$
\end{definition}

\begin{definition}
We can also define a related invariant, associated this time to just
one given point $x\in X$:
$$p(L,x) := {\rm max}\{l~|~L\otimes m_x^l~{\rm is ~}M{\rm -regular}\}.$$
The definition does not depend on $x$ because of the homogeneity of $X$, and
so we will denote this invariant simply by $p(L)$.
\end{definition}

Our main interest will be in the asymptotic versions of these indices, which
turn out to be related to the Seshadri constant associated to $L$.

\begin{definition}
The \emph{asymptotic $M$-regularity index} of $L$ and its punctual counterpart
are defined as
$$\rho(L) := \underset{n}{{\rm sup}} \frac{m(L^n)}{n}$$
and
$$\rho^\prime (L) := \underset{n}{{\rm sup}} \frac{p(L^n)}{n}.$$
\end{definition}

\noindent
The main result of this section is:

\begin{theorem}\label{seshadri}
We have the following inequalities:
$$\epsilon(L) = \rho^{\prime}(L) \geq \rho(L)\geq 1.$$
In particular $\epsilon(L)\ge {\rm max}\{m(L),1\}$.
\end{theorem}

This improves a result of Nakamaye (cf. \cite{nakamaye} and the references
therein). Nakamaye also shows that
$\epsilon(L) = 1$ for some line bundle $L$ if and only if $X$ is
the product of an elliptic curve with another abelian variety, so
then a similar result holds for the invariant $\rho^{\prime}(L)$.
As explained in \cite{pp2} \S3, the value of
$m(L)$ is reflected in the geometry of the map to projective space
given by $L$. Here is a basic example:

\begin{example}
If $L$ is very ample -- or more generally gives a birational morphism
outside a codimension $2$ subset --  then $m(L)\ge 2$, and so by the
theorem above $\epsilon (L)\ge 2$. Note that on an arbitrary smooth
projective variety very ampleness implies in general only that $\epsilon(L,x)\ge 1$
at each point.
\end{example}

The proof of Theorem \ref{seshadri} is a simple application of the
$M$-regularity criterion \ref{M_reg}, via the results of \cite{pp2} \S3.
To understand the growth of the usual invariants we use the
relationship with the notions of $k$-jet ampleness and separation of jets. Namely let's
denote by $s(L,x)$ the largest number $s\ge 0$ such that $L$ separates
$s$-jets at $x$. Recall also the following:

\begin{definition}
A line bundle $L$ is called $k$-\emph{jet ample}, $k\geq 0$, if the restriction map
$$H^0(L)\longrightarrow H^0(L\otimes \OO_X/
m_{x_1}^{k_1}\otimes \ldots \otimes m_{x_p}^{k_p})$$
is surjective for any distinct points $x_1,\ldots,x_p$ on $X$ such that $\Sigma k_i =k+1$.
\end{definition}

Note that if $L$ is $k$-jet ample, then it separates
$k$-jets at every point. Recall from \cite{pp2} Theorem 3.8 and
Proposition 3.5 the following facts:

\begin{proposition}\label{both_ways}
(i) $L^n$ is ($n + m(L) -2$)-jet ample, so in particular $s(L^n,x)\ge n + m(L) -2$.
\newline
\noindent
(ii) If $L$ is $k$-jet ample, then $m(L)\geq k+1$.
\end{proposition}

This immediately points in the direction of local positivity, since
one way to interpret the Seshadri constant of $L$ is (independently of $x$):
$$\epsilon (L) = \underset{n}{{\rm limsup}} \frac{s(L^n,x)}{n} =
\underset{n}{{\rm sup}} \frac{s(L^n,x)}{n}.$$
The last equality follows from the fact that jet-separation satisfies
the well-known superadditivity relation $s(L_1\otimes L_2,x)\geq
s(L_1,x) + s(L_2,x)$ for any two line bundles $L_1$ and $L_2$ on $X$.
To establish the connection with the asymptotic invariants above
we also need the following

\begin{lemma}\label{comparison}
For any $n\geq 1$ and any $x\in X$ we have $s(L^{n+1},x)\geq m(L^n)$.
\end{lemma}
\begin{proof}
This follows immediately from the $M$-regularity criterion Theorem \ref{M_reg}: if
$L^n\otimes m_{x_1}^{k_1}\otimes \ldots \otimes m_{x_p}^{k_p}$ is
$M$-regular, then $L^{n+1}\otimes m_{x_1}^{k_1}\otimes \ldots
\otimes m_{x_p}^{k_p}$ is globally generated, and so by \cite{pp2} Lemma 3.3,
$L^{n+1}$ is $m(L)$-jet ample.
\end{proof}

\begin{proof}(\emph{of Theorem \ref{seshadri}.})
Note first that for every $p\ge 1$ we have
\begin{equation}\label{scaling}
m(L^n)\ge m(L) + n -1,
\end{equation}
which follows immediately from the two parts of Proposition \ref{both_ways}.
In particular $m(L^n)$ is always at least $n-1$, and so $\rho(L)\ge 1$.

Putting together the definitions, (\ref{scaling}) and
Lemma \ref{comparison}, we obtain the main inequality $\epsilon(L)
\geq \rho(L)$.
Finally, and less surprisingly, we see equally quickly that the
asymptotic punctual index computes precisely the Seshadri constant.
Indeed, by completely similar arguments as above, we have that
for any ample line bundle $L$ and any $p\geq 1$ one has
$$p(L^n)\ge s(L^n,x)~{\rm and~} s(L^{n+1},x)\ge p(L^n, x).$$
The statement follows then from the definition.
\end{proof}

\begin{remark}
What the proof above shows is that one can give an interpretation
for $\rho(L)$ similar
to that for $\epsilon (L)$ in terms of separation of jets.
In fact $\rho(L)$ is precisely the "asymptotic jet ampleness" of $L$ (stronger
then jet separation), namely:
$$\rho (L) = \underset{n}{{\rm sup}} \frac{a(L^n)}{n},$$
where $a(M)$ is the largest
integer $k$ for which the line bundle $M$ is $k$-jet ample.
\end{remark}

In this respect, an interesting -- though admittedly quite optimistic --
question is whether the asymptotic $M$-regularity index computes precisely the Seshadri constant:

\begin{question}
Do we always have $\epsilon(L) = \rho(L)$?
\end{question}

\begin{remark}
Another, rather surprising, lower bound for the Seshadri constant
of a polarization on an abelian variety has been given by Lazarsfeld in \cite{rob1}.
This is expressed in terms of a metric invariant defined by Buser and
Sarnak. In comparison, Lazarsfeld's result can be made effective
based on the Buser-Sarnak result giving a lower bound for that
particular invariant (cf. \emph{loc. cit.}). For the same reason,
it would be of
considerable interest to find an independent lower bound for the
asymptotic invariant $\rho(L)$.
\end{remark}

\begin{question}
Can one give independent lower bounds for $\rho(L)$ or $\rho^\prime(L)$?
\end{question}

It may be possible to do this for generic abelian varieties by
constructing specific examples.
In the other direction, there are numerous bounds on Sesahdri
constants, which in turn give bounds for the $M$-regularity
indices that (at least to us) are not obvious from the definition.
Essentially each of the results listed in \cite{rob2} Ch.I \S5
gives some sort of bound. Let's just give a couple of examples:

\begin{corollary}
If $(J(C), \Theta)$ is a Jacobian with the usual principal
polarization, then $m(n\Theta)\le \sqrt{g}\cdot n$. On an
arbitrary abelian variety, for any principal polarization
$\Theta$ we have $m(n\Theta)\le (g!)^{\frac{1}{g}}\cdot n$.
\end{corollary}
\begin{proof}
It is shown in \cite{rob1} that $\epsilon(\Theta)\le \sqrt{g}$.
We then apply Theorem \ref{seshadri}. For the other
bound we use the usual elementary upper bound for Seshadri constants,
namely $\epsilon(\Theta) \leq (g!)^{\frac{1}{g}}$.
\end{proof}

\begin{corollary}
If $(A, \Theta)$ is a very general PPAV, then there exists
at least one $n$ such that $p(n\Theta)\ge \frac{2^{\frac{1}{g}}}{4}
(g!)^{\frac{1}{g}}\cdot n$.
\end{corollary}
\begin{proof}
Here we use the lower bound given in \cite{rob1} via the
Buser-Sarnak result.
\end{proof}

There are of course more specific results on $\epsilon (\Theta)$ for Jacobians (cf. \cite{debarre}
Theorem 7), each giving a corresponding result for $m(n\Theta)$.
It would be more satisfactory to have a concrete answer in this case, but note that in this generality
the problem should be quite difficult, since it would also answer conjectures about the Seshadri
constant (for example the fact that $\epsilon(\Theta) < 2$ chracterizes hyperelliptic Jacobians).

\begin{question}
Can we understand $m(n\Theta)$ individually on Jacobians, at least for small $n$, in
terms of the geometry of the curve?
\end{question}

As a simple example, the question above has a clear answer for elliptic curves.
We know that on an elliptic curve $E$ a line bundle $L$
is $M$-regular if and only if ${\rm deg}(L)\ge 1$, i.e. if
and only if $L$ is ample. From the definition of $M$-regularity
we see then that if ${\rm deg}(L)=d>0$, then $m(L)=d-1$.
This implies that on an elliptic curve $m(n\Theta) = n-1$ for all $n\ge 1$.
However, this is misleading when we look at the case of curves of higher genus. In fact
in the simplest case we have the following general:

\begin{proposition}
If $(X,\Theta)$ is an irreducible principally polarized abelian variety of dimension at least $2$, then $m(2\Theta)\geq 2$.
\end{proposition}
\begin{proof}
This is an immediate consequence of the existence of the Kummer map. The linear series $|2\Theta|$
induces a $2:1$ map of $X$ onto its image in $\PP^{2^g-1}$, with injective differential. Thus the
cohomological support locus for $\OO(2\Theta)\otimes m_x\otimes m_y$ consists of two points, while
the one for $\OO(2\Theta)\otimes m_x^2$ is empty.
\end{proof}

\section{Regularity of Picard bundles and vanishing on symmetric products}

In this section we study the regularity of Picard bundles over the
Jacobian of a curve, in other words we give a quantitative estimate for their
positivity with respect to the natural polarization. Our study is not a direct
consequence of $M$-regularity, but integrates nicely in the context of its strong
version called Theta regularity (cf. \S2).
The point is that one can prove vanishing results by combining the Fourier-Mukai
transform with classical resolution type methods for determinantal
varieties, involving in particular the Eagon-Northcott complex.

Let $C$ be a smooth curve of genus $g\geq 2$, and denote by $J(C)$ the
Jacobian of $C$,  and by $C_n$ the $n$-th symmetric product of $C$.
The objects we are interested in are the Picard
bundles on $J(C)$: for a given $n\geq 2g-1$, a Picard bundle is loosely speaking
a vector bundle over ${\rm Pic}^n(C)$
whose projectivization is $C_n$, so that the projectivizations
of its fibers parametrize all linear series of degree $n$ on $C$ (cf. \cite{acgh} Ch.VII \S2).
We will look at such a bundle $E$ on $J(C)$, via translating by a line bundle $L\in {\rm Pic}^n(C)$,
so we make the convention that the vector bundle fiber of $E$ at a point $\xi\in J(C)$
is $H^0(L\otimes \xi)$. We will somewhat abusively call this the \emph{$n$-th Picard
bundle of $C$}.

\begin{proposition}\label{fm_picard}
For any $k\geq 1$, let $\pi_k: C^k \rightarrow J(C)$ a desymmetrized
Abel-Jacobi mapping and let $L$ be a line bundle on $C$ of degree $n>>0$ as above.
Then ${\pi_k}_*(L\boxtimes \ldots\boxtimes L)$ satisfies I.T. with index $0$, and
$$({\pi_k}_*(L\boxtimes \ldots\boxtimes L))^{\widehat{}} = \otimes^k E,$$
where $E$ is the $n$-th Picard bundle of $C$.
\end{proposition}
\begin{proof}
This is a generalization of the well-known fact (cf. \cite{mukai1} \S4) that the
Picard bundle $E$ is the Fourier transform of $i_* L$, where $i$ is an Abel-Jacobi
embedding of $C$ in the Jacobian, and we only briefly sketch the proof. Indeed, the positivity
of $L$ insures the fact that ${\pi_k}_*(L\boxtimes \ldots\boxtimes L)$ satisfies
I.T. with index $0$, as the fibers of the second projection to the dual Jacobian
have no higher cohomology (this can be easily seen using the K\" unneth formula).
This implies that the fibers of the Fourier transform are naturally isomorphic
to $H^0(C,L\otimes\xi)^{\otimes k}$, which characterizes $\otimes^k E$.
\end{proof}

The following theorem is the main cohomological result we are aiming for.
The point to keep in mind is that Picard bundles are known to be
negative (i.e with ample dual bundle), so vanishing theorems are not automatic.
We give an effective range for achieving the strongest vanishing one can hope for
with the smallest possible "positive" perturbation. To be very precise, everything
that follows holds if $n$ is assumed to be at least $4g-4$. (However the value of $n$
does not affect the applications.)

\begin{theorem}\label{mreg_picard}
For every $1\leq k\leq g-1$, $\otimes^k E$ is strongly $2$-$\Theta$-regular.
\end{theorem}
\begin{proof}\footnote{We
are grateful to Olivier Debarre for pointing out a numerical mistake in the statement, in a
previous version of this paper.}
We will use loosely the notation $\Theta$ for any translate of the canonical theta divisor.
The statement of the theorem becomes then equivalent to the vanishing
$$h^i(\otimes^k E\otimes \OO(\Theta))=0,~\forall~ i>0, ~\forall~ 1\leq k\leq g-1.$$
To prove this vanishing we use the Fourier-Mukai transform. The first point is
that Proposition \ref{fm_picard} above, combined
with Mukai's main duality theorem \cite{mukai1} Theorem 2.2, tells us precisely
that  $\otimes^k E$ satisfies W.I.T. with index $g$, and its Fourier transform is
$$\widehat{\otimes^k E}=(-1_J)^*{\pi_k}_* (L\boxtimes\ldots\boxtimes L).$$

The next point is that the cohomology groups involved
can be computed on the dual Jacobian via the Fourier transform. We have the following
sequence of isomorphisms:
$$H^i(\otimes^k E\otimes \OO(\Theta))\cong {\rm Ext}^i(\OO(-\Theta), \otimes^k E)
\cong {\rm Ext}^i(\widehat{\OO(-\Theta)}, \widehat{\otimes^k E})$$
$$\cong {\rm Ext}^i(\OO(\Theta), (-1_J)^*{\pi_k}_*(L\boxtimes\ldots\boxtimes L))\cong
H^i((-1_J)^*{\pi_k}_*(L\boxtimes\ldots\boxtimes L)\otimes \OO(-\Theta)).$$
Here we are using the correspondence between the ${\rm Ext}$ groups given in
\cite{mukai1} Corollary 2.5, plus the fact that both $\OO(-\Theta)$ and $\otimes^k E$ satisfy
W.I.T. with index $g$ and that $\widehat{\OO(-\Theta)}=\OO(\Theta)$.

As we are loosely writing $\Theta$ for any translate, multiplication by $-1$ does
not influence the vanishing, so the result follows if we show:
$$h^i({\pi_k}_*(L\boxtimes\ldots\boxtimes L)\otimes \OO(-\Theta))=0, ~\forall~i>0.$$
Now the image $W_k$  of the Abel-Jacobi map $u_k:C_k \rightarrow J(C)$ has rational
singularities (cf. \cite{kempf2}), so we only need to prove the vanishing:
$$h^i(u_k^*({\pi_k}_*(L\boxtimes\ldots\boxtimes L)\otimes \OO(-\Theta)))=0, ~\forall~i>0.$$
Thus we are interested in the skew-symmetric part of the cohomology group $H^i(C^k, (L\boxtimes\ldots\boxtimes L)\otimes \pi_k^*\OO(-\Theta))$, or, by Serre duality that of
$$H^i(C^k, ((\omega_C\otimes L^{-1})\boxtimes\ldots\boxtimes (\omega_C\otimes L^{-1}))
\otimes \pi_k^*\OO(\Theta)), ~{\rm for}~i<k.$$
At this stage we can essentially invoke a Serre vanishing type argument, but it is
worth noting that the computation can be in fact made very concrete. For the identifications
used next we refer to \cite{izadi} Appendix 3.1. As $k\le g-1$, we can write
$$\pi_k^*\OO(\Theta)\cong ((\omega_C\otimes A^{-1})\boxtimes\ldots\boxtimes
(\omega_C\otimes A^{-1}))\otimes \OO(-\Delta),$$
where $\Delta$ is the union of all the diagonal divisors in $C^k$ and $A$ is a line
bundle of degree $g-k-1$. Then the skew-symmetric part of the cohomology groups we are looking at is isomorphic to
$$S^i H^1(C, \omega_C^{\otimes 2}\otimes A^{-1}\otimes L^{-1})  \otimes \wedge^{k-i}
H^0(C, \omega_C^{\otimes 2}\otimes A^{-1}\otimes L^{-1}),$$
and since for $1\leq k\leq g-1$ and $n\ge 4g-4$ the degree of the line bundle  $\omega_C^{\otimes 2}\otimes A^{-1}\otimes L^{-1}$ is negative, this vanishes precisely for $i<k$.
\end{proof}

\bigskip

An interesting consequence of the vanishing result for Picard bundles proved above
is a new -- and in some sense more classical -- way to deduce Theorem 4.1 of
\cite{pp1} on the regularity of the ideal sheaves $\I_{W_d}$ on the Jacobian $J(C)$.
This theorem has a number of nice consequences on the equations of the $W_d$'s
-- in particular on those of the curve $C$ -- inside $J(C)$, and also to some useful
vanishing results for pull-backs  of theta divisors to symmetric products. For this circle
of ideas we refer the reader to \cite{pp1} \S4.

For any $1\leq d\leq g-1$, $g \geq 3$, consider $u_d :C_d\longrightarrow J(C)$
to be an Abel-Jacobi mapping of the symmetric product (depending on the choice
of a line bundle of degree $d$ on $C$), and denote by $W_d$
the image of $u_d$ in $J(C)$.

\begin{theorem}\label{two_theta}
For every $1\leq d\leq g-1$, the ideal sheaf $\I_{W_d}$ is strongly $3$-$\Theta$-regular.
\end{theorem}
\begin{proof}
We have to prove that:
$$h^i(\I_{W_d}\otimes \OO(2\Theta)\otimes \alpha) = 0,~\forall~i>0,~\forall~\alpha\in
{\rm Pic}^0(J(C)).$$
In the rest of the proof, by $\Theta$ we will understand generically
any translate of the canonical theta divisor, and so $\alpha$ will disappear from the notation.

It is well known that $W_d$ has a natural determinantal structure, and its ideal
is resolved by an Eagon-Northcott complex. We will chase the vanishing along this
complex. This setup is precisely the one used by Fulton and Lazarsfeld in order to prove
for example the existence theorem
in Brill-Noether theory, and for explicit details on this circle of ideas
we refer to \cite{acgh} Ch.VII \S2.

Concretely, $W_d$ is the "highest" degeneracy locus of a map of vector bundles
$$\gamma: E\longrightarrow F,$$
where ${\rm rk}F=m$ and ${\rm rk}E=n=m+d-g+1$, with $m>>0$ arbitrary. The bundles
$E$ and $F$ are well understood: $E$ is the $n$-th Picard bundle of $C$,
discussed above, and $F$ is a direct
sum of topologically trivial line bundles. (For simplicity we are again moving the
whole construction on $J(C)$ via the choice of a line bundle of degree $n$.) In other
words, $W_d$ is scheme theoretically the locus where the dual map
$$\gamma^*:F^*\longrightarrow E^*$$
fails to be surjective. This locus is resolved by an Eagon-Northcott complex
(cf. \cite{kempf}) of the form:
$$0\rightarrow \wedge^m F^*\otimes S^{m-n}E\otimes {\rm det}E\rightarrow \ldots
\rightarrow \wedge^{n+1} F^*\otimes E\otimes {\rm det}E\rightarrow \wedge^n F^*
\rightarrow \I_{W_d}\rightarrow 0.$$
As it is known that the determinant of $E$ is nothing but $\OO(-\Theta)$, and
since $F$ is a direct sum of topologically trivial line bundles, the statement of
the theorem follows by chopping this into short exact sequences, as long as we prove:
$$h^i(S^k E\otimes \OO(\Theta))=0,~\forall~ i>0, ~\forall~ 1\leq k\leq m-n=g-d-1.$$
Since we are in characteristic zero, $S^k E$ is naturally a direct summand in $\otimes^k E$,
and so it is sufficient to prove that:
$$h^i(\otimes^k E\otimes \OO(\Theta))=0,~\forall~ i>0, ~\forall~ 1\leq k\leq g-d-1.$$
But this follows from Theorem \ref{mreg_picard}.
\end{proof}

\begin{remark}
Using \cite{pp1} Proposition 2.9, we have a strong version of the "abelian" Castelnuovo-Mumford
Lemma \cite{pp1} Theorem 6.3, namely a strongly $m$-$\Theta$-regular sheaf is strongly $k$-$\Theta$-regular for every $k\geq m$. Thus Theorem \ref{mreg_picard} and Theorem \ref{two_theta} imply
that $\otimes^k E$ is strongly $k$-$\Theta$-regular for every $k\geq 2$ and $\I_{W_d}$ is strongly $k$-$\Theta$-regular for every $k\geq 3$.
\end{remark}

\begin{question}
An interesting question, extending the result above, is the following: what is the $\Theta$-regularity of
the ideal of an arbitrary Brill-Noether locus $W_d^r$?
\end{question}

We describe below one case in which the answer can already be
given, namely that of the singular locus of the Riemann theta
divisor on a non-hyperelliptic jacobian. It should be noted that
in this case we do not have strong $3$-$\Theta$-regularity any
more (but rather only strong $4$-$\Theta$-regularity, by the same
\cite{pp1} Proposition 2.9).

\begin{proposition}\label{singtheta}
The ideal sheaf ${\mathcal I}_{W^1_{g-1}}$ is 3-$\Theta$-regular.
\end{proposition}
\begin{proof} In fact it follows from the results of \cite{vgi}
that $$ h^i({\mathcal I}_{W^1_{g-1}}\otimes {\mathcal
O}(2\Theta)\otimes\alpha)=\begin{cases}0&\hbox{for $i\ge g-2$,
$\forall\alpha\in {\rm Pic}^0(J(C))$}\cr 0&\hbox{for $0<i<g-2$,
$\forall\alpha\in {\rm Pic}^0(J(C))$ such that $\alpha\ne
{\mathcal O}_{J(C)}.$}\cr
\end{cases}$$
For the reader's convenience, let us briefly recall the
relevant points from Section 7 of \cite{vgi}. We denote for simplicity, via translation,
$\Theta=W_{g-1}$, (so that $W^1_{g-1}= {\rm Sing}(\Theta$)). In the first
place, from the exact sequence
$$0\rightarrow {\mathcal O}(2\Theta)\otimes\alpha\otimes \OO(-\Theta)\rightarrow {\mathcal
I}_{W^1_{g-1}}(2\Theta)\otimes\alpha\rightarrow {\mathcal
I}_{W^1_{g-1}/\Theta}(2\Theta)\otimes\alpha\rightarrow 0$$ it
follows that
$$h^i( J(C) ,{\mathcal
I}_{W^1_{g-1}}(2\Theta)\otimes\alpha)=h^i(\Theta, {\mathcal
I}_{W^1_{g-1}/\Theta}(2\Theta)\otimes\alpha)\hbox{ for $i>0$.}$$
Hence one is reduced to a computation on $\Theta$.  It is a
standard fact  (see e.g. \cite{vgi}, 7.2) that, via the
Abel-Jacobi map $u=u_{g-1}:C_{g-1}\rightarrow \Theta\subset J(C)$,
$$h^i(\Theta, {\mathcal
I}_{W^1_{g-1}/\Theta}(2\Theta)\otimes\alpha)=h^i(C_{g-1},
L^{\otimes 2}\otimes \beta\otimes {\mathcal I}_Z),$$ where $Z=u^{-1}(W^1_{g-1})$,
$L=u^*{\mathcal O}_X(\Theta)$ and $\beta=u^*\alpha$. We now use
the standard exact sequence (\cite{acgh}, p.258):
$$0\rightarrow T_{C_{g-1}}\buildrel{du}\over\rightarrow
H^1(C,{\mathcal O}_C)\otimes{\mathcal O}_{C_{g-1}}\rightarrow
L\otimes {\mathcal I}_Z\rightarrow 0.$$
Tensoring with $L\otimes
\beta$, we see that it is sufficient to prove that
$$H^i(C_{g-1},T_{C_{g-1}}\otimes L\otimes\beta )=0, \> \forall i\ge
2, \> \forall \beta\ne {\mathcal O}_{C_{g-1}}.$$
To this end we use the well known fact (cf. \emph{loc. cit.}) that
$$T_{C_{g-1}}\cong p_*{\mathcal O}_D(D)$$
where $D\subset C_{g-1}\times C$ is the universal divisor and $p$
is the projection onto the first factor. As $p_{|D}$ is finite, the
degeneration of the Leray spectral sequence and the projection formula
ensure that
$$h^i(C_{g-1}, T_{C_{g-1}}\otimes L\otimes\beta)=h^i(D, {\mathcal
O}_D(D)\otimes p^*(L\otimes \beta)),$$
which are zero for $i\ge 2$ and $\beta$ non-trivial by \cite{vgi}, Lemma 7.24.
\end{proof}

\section{Pluricanonical maps of irregular varieties of maximal
Albanese dimension}

It is well known that the minimal pluricanonical embedding working
for all smooth curves of general type is the tricanonical one. It
turns out that this very classical result has a generalization to
arbitrary dimension. In fact Chen and Hacon (\cite{chenhacon},
Theorem 4.4)  recently proved that the tricanonical map of a
smooth complex irregular variety of general type $Y$, having
generically finite Albanese map and such that $\chi(\omega_Y)>0$,
is birational onto its image. The main point of this section is
that the concept of M-regularity, combined with well-known results
of Ein, Green and Lazarsfeld, provides -- under mildly more
restrictive hypotheses -- a very quick and conceptually simple
proof of a more explicit version of this statement, as well as of
other related facts. To put things into perspective, let us recall
that, by a theorem of Ein, Lazarsfeld and Green (\cite{einlaz},
Theorem 3), partly conjectured by Koll\'ar, given a smooth variety
of general type $Y$ with generically finite Albanese map,  then
$\chi(\omega_Y)\ge 0$, and if equality holds then the Albanese
image of $Y$ is ruled by tori. We show the following:

\begin{theorem}\label{3can}
Let $Y$ be a smooth projective complex minimal variety of general type
such that  its Albanese map $a:Y\rightarrow {\rm Alb}(Y)$ is generically
finite. If $a(Y)$ is not ruled by tori, then $\omega_Y^{\otimes
3}$  is very ample on the open set  $a^{-1}(T)$, where $T$ is the
open set of points of $a(Y)$ over which the fiber of $a$ is  finite.
\end{theorem}

In preparation for the proof, let us settle some preliminary results.
First we introduce a slight generalization of Definition \ref{continuous}.

\begin{definition}
Let $a:Y\rightarrow X$ be a map from a projective variety $Y$ to an
abelian variety $X$ and let ${\mathcal F}$ be a coherent
sheaf on $Y$. Let also $Z\subset Y$ be a closed subset. We say that
${\mathcal F}$ is \emph{continuously globally generated away from
$Z$ (with respect to the map $a$)} if, for any Zariski open set
$U\subset \Pic0$, $Z$ does not meet  the support of the cokernel of
the evaluation map
$${\rm ev}_U:\bigoplus_{\alpha\in U}H^0({\mathcal F}\otimes
a^*\alpha)\otimes a^*\alpha^\vee\rightarrow {\mathcal F}.$$
\end{definition}
With this terminology,  we have the following trivial
generalization of Theorem \ref{M_reg}(b) (the proof is the same).

\begin{proposition}
 Let ${\mathcal F}$ and ${\mathcal L}$ be
respectively a coherent sheaf and an invertible sheaf on $Y$. If they
are both continuously globally generated  away from $Z$ (with
respect to $a$) then ${\mathcal F}\otimes L\otimes a^*\alpha$ is
globally generated away from $Z$ for any $\alpha\in\Pic0$.
\end{proposition}

From this point on, $Y$ will be a smooth projective variety,
the map $a$ will be the Albanese map of $Y$, and
the subset $Z$ will be the inverse image, via $a$, of the locus of
points of $a(Y)$ having non-finite fiber. The key point is

\begin{lemma}\label{key_canonical}
If $\dim Y=\dim a(Y)$ and $a(Y)$ is not ruled by tori then:
\newline\noindent (i) $a_*\omega_Y$ is an M-regular sheaf on $X={\rm Alb}(Y)$;
\newline\noindent (ii) $a_*\omega_Y$ is continuously globally
generated;
\newline\noindent (iii) $\omega_Y$ is continuously globally generated away from
$Z$;
\newline\noindent (iv) $\omega_Y^{\otimes 2}\otimes a^*\alpha$ is globally generated
away from $Z$ for any $\alpha\in \Pic0$.
\end{lemma}
\begin{proof}
(i) Given a coherent sheaf ${\mathcal F}$ on $Y$, let us denote
$$V_i({\mathcal F})=\{\alpha\in \Pic0 \>|\>
h^i({\mathcal F}\otimes a^*\alpha)>0\},$$
the $i$-th \emph{cohomological support locus} of $\F$.
As we are assuming that $a$
is generically finite onto its image, by the Generic Vanishing
Theorem of \cite{greenlaz} we have that ${\rm codim}~ V_i({\omega_Y})\ge i$ for all
$i$. Moreover, by a well-known argument of \cite{einlaz} (end of the proof
of Theorem 3), if $a(X)$ is not ruled by tori then ${\rm codim}~
V_i(\omega_Y)>i$ for all $i\ge 1$. Since, by
Grauert-Riemenschneider vanishing, $R^ia_*\omega_Y=0$ for any
$i\ne 0$, the projection formula gives
$$H^i(Y,\omega_Y\otimes a^*\alpha)\cong H^i(X,(a_*\omega_Y)\otimes
\alpha).$$
This implies that ${\rm codim}~ V_i(a_*\omega_X)> i$ for all
$i>0$, which implies by base change that $a_*\omega_X$ is M-regular.

\noindent
(ii) This follows from (i) via the $M$-regularity criterion Theorem \ref{M_reg}.

\noindent
(iii) It is immediate to see that, as with global generation, continuous global
generation is preserved by finite maps. We then apply (ii).

\noindent
(iv) This follows from (iii) and the previous Proposition.
\end{proof}

\begin{remark}
The concept of $M$-regularity in the case of canonical bundles on varieties of maximal
Albanese dimension interprets very
conceptually some of the results of \cite{einlaz}. It is in fact easy to see, as noted in
\cite{pp1}, that a non-zero  $M$-regular sheaf has global sections. Thus by
the definition, if $a_*\omega_Y$ is $M$-regular, then $\chi (\omega_Y)>0$. So Lemma
\ref{key_canonical} reinterprets the result of Ein and Lazarsfeld mentioned above, using
of course half of their argument, employed in the proof of \ref{key_canonical} (i).

\begin{corollary}(\cite{einlaz}, Theorem 3)
If $Y$ is a smooth projective variety of maximal Albanese dimension with $\chi(\omega_Y)=0$,
then the Albanese image $a(Y)$ is ruled by tori (since $\omega_Y$ is not $M$-regular). 
\end{corollary}

\end{remark}

\begin{proof} (\emph{of Theorem \ref{3can}})
The statement is equivalent to the fact that, for any $y\in Y-Z$,
the sheaf $I_y\otimes \omega_Y^{\otimes 3}$ is globally generated
away from $Z$. Consider such a $y$. As in deducing (iv) in the previous
Lemma, via the $M$-regularity criterion the
statement will follow from the fact that the
sheaf $a_*(I_y\otimes \omega_Y^{\otimes 2})$ is M-regular,
which is what we will prove.

To this end, note that by Lemma \ref{key_canonical} (iv),
$\omega_Y^2\otimes a^*\alpha$ is globally generated at $y$.
Since, by Kawamata-Viehweg vanishing,
$H^i(Y,\omega_Y^2\otimes a^*\alpha)=0$ for any $i>0$ (here we are
using the hypothesis that  $\omega_Y$ is big and nef), it follows
easily that
$$H^i(I_y\otimes\omega_Y^{\otimes 2} \otimes
a^*\alpha)=0,~{\rm for~ any~} \alpha\in \Pic0.$$
On the other hand we have
that $R^ia_*(I_y\otimes \omega_Y^2)=0$ for any $i>0$. This
follows from the standard exact sequence
$$0\rightarrow
I_y\otimes\omega^{\otimes 2}_Y\rightarrow\omega^{\otimes
2}_Y\rightarrow {\omega_Y^{\otimes 2}}_{|y}\rightarrow 0,$$
using the fact that, by a well-known generalization of
Grauert-Riemenschneider vanishing (\cite{kollarmori}, Corollary 2.8),
$R^ia_*(\omega_Y^{\otimes 2})=0$ for any $i>0$. Therefore, by the
projection formula, $a_*(I_y\otimes\omega_Y^{\otimes 2})$
satisfies I.T. with index $0$, so in particular it is M-regular. This proves the
Theorem.
\end{proof}

The following generalization of Theorem \ref{3can} to higher order jets is proved 
the same way (cf. also \cite{pp2} \S3). 

\begin{theorem}
Let $Y$ be a smooth projective complex minimal variety of general
type such that  its Albanese map $a:Y\rightarrow {\rm Alb}(Y)$ is
generically finite. If $a(Y)$ is not ruled by tori, then, for any
$k\ge 0$ , $\omega_Y^{\otimes k+2}$  is k-jet ample on the open
set $a^{-1}(T)$, where $T$ is the open set of points of $a(Y)$
over which the fiber of $a$ is  finite.
\end{theorem}

To complete the picture, we recall the following general
result, indicated in \cite{pp1}, Remark 5.3, and partly obtained
independently by Chen and Hacon (compare e.g. \cite{chenhacon},
Corollaries 3.3 and 4.3). We omit the proof, which is very
similar to the previous one.

\begin{theorem}
Let $Y$ be a smooth projective complex  variety such that  its
Albanese map $a:Y\rightarrow {\rm Alb}(Y)$ is generically finite,
and let $L$ be a big and nef line bundle on $Y$. Then, for any
$k\ge 0$, $(\omega_Y\otimes L)^{\otimes k+2}$ is k-jet ample on
the open set $a^{-1}(T)$, where $T$ is the open set of points of
$a(Y)$ over which the fiber of $a$ is finite. For instance, if $Y$
is minimal of general type, then $\omega_Y^{\otimes 6}$ is very
ample on the open set $a^{-1}(T)$.
\end{theorem}

\section{Numerical study of semihomogeneous vector bundles}

An idea that originated in work of Mukai is that on abelian varieties the
class of vector bundles to which the theory of line bundles should generalize
naturally is that of \emph{semihomogeneous} bundles (cf. \cite{mukai1},
\cite{mukai3}, \cite{mukai4}). These vector bundles are
semistable, behave nicely under isogenies and Fourier transforms, and have a
Mumford type theta group theory as in the case of line bundles (cf. \cite{umemura}).
The purpose of this section is to show that this analogy can be extended to include
effective global generation and normal generation statements dictated by specific numerical
invariants measuring positivity.

The guiding problem is the
following: given a semihomogeneous bundle $E$ on an abelian variety $X$,
such that $h^0(E)\neq 0$,
find sufficient "positivity" conditions on $E$ which ensure that $E$ is globally or
normally generated.  Recall that the later is Mumford's terminology for the surjectivity of the
multiplication map $H^0(E) \otimes H^0(E)\rightarrow H^0(E^{\otimes 2})$ -- the model to keep in mind is the global generation of $\OO(2\Theta)$ and the projective normality of $\OO_X(3\Theta)$ for any ample divisor $\Theta$ on $X$.
We would like to integrate this into the general regularity
theory of \cite{pp1}, since $\Theta$-regularity as described in \S2 is precisely a measure of the
(cohomological) positivity of $E$. As a minimal requirement we have to
ask that $E$ be $0$-$\Theta$-regular (recall that this means that
$E(-\Theta)$ is $M$-regular). Note that this automatically implies $h^0(E)\neq 0$.

\noindent
{\bf Basics on semihomogeneous bundles.}
Let $X$ be an abelian variety of dimension $g$ over an algebraically
closed field. As a general convention, for a numerical class $\alpha$ we will use the
notation $\alpha >0$ to express the fact that $\alpha$ is ample. If the class
is represented by an effective divisor, then the condition of being ample
is equivalent to $\alpha^g >0$.

\noindent
Let $L$ be a line bundle on $X$. We denote by $\phi_L$
the isogeny defined by $L$:
$$\begin{array}{cccc}
\phi_{L}:& X& \longrightarrow & {\rm Pic}^{0}(X)\cong\widehat{X}\\
&x &\leadsto& t_{x}^*L\otimes L^{-1}.
\end{array}
$$
Note that if $L$ is a principal polarization, $\phi_L$ is just a
self-duality of $X$.

\smallskip

\begin{definition}(\cite{mukai3})
A vector bundle $E$ on $X$ is called \emph{semihomogeneous} if for every $x\in X$,
$t_x^*E\cong E\otimes \alpha$, for some $\alpha\in {\rm Pic}^0(X)$.
\end{definition}

It is a general principle, described later in this subsection,
that the study of arbitrary semihomogeneous bundles can be
reduced to that of \emph{simple} semihomogeneous ones, i.e those with no
nontrivial endomorphisms.
We will use a few basic properties of simple semihomogeneous bundles which can
be found in \cite{mukai3}.
In the next lemmas $E$ is a simple semihomogeneous of rank $r$ on $X$.

\begin{lemma}\label{sh_pullback}(\cite{mukai3} Proposition 7.3)
There exists an isogeny $\pi: Y\rightarrow X$ and a line bundle $M$ on $Y$
such that
$$\pi^*E\cong \underset{r}{\bigoplus}M.$$
\end{lemma}

\begin{lemma}\label{sh_pushforward}(\cite{mukai3} Theorem 5.8(iv))
There exists an isogeny $\phi: Z\rightarrow X$ and a line bundle $L$ on $Z$
such that
$$\phi_* L=E.$$
\end{lemma}

The second lemma implies that any simple semihomogeneous bundle satisfies
an Index Theorem analogous to the line bundle one (cf. \cite{mumford1} \S16).

\begin{lemma}\label{index_theorem}
Let $E$ be a nondegenerate (i.e. $\chi(E)\neq 0$) simple semihomogeneous bundle
on $X$. Then exactly one cohomology group $H^i(E)$ is nonzero.
\end{lemma}
\begin{proof}
This follows immediately from the similar property of the line bundle $L$
in Lemma \ref{sh_pushforward}.
\end{proof}

We will be concerned with semihomogeneous bundles which have some sort
of positivity, so in particular are nondegenerate and have global
sections. A first property is that they have well-defined Fourier-Mukai transforms.

\begin{lemma}\label{ss_it0}
Assume $E$ is nondegenerate simple semihomogeneous, such that $h^0(E)\neq 0$.
Then $E$ satisfies I.T. with index $0$.
\end{lemma}
\begin{proof}
For every $\alpha \in \Pic0$, there exists an $x\in X$ such that $E\otimes \alpha
\cong t_x^* E$. This implies that
$$H^0(E\otimes \alpha) \cong H^0(t_x^*E)\cong H^0(E)\neq 0.$$
By Lemma \ref{index_theorem} this means that $h^i(E\otimes \alpha) =0$ for all
$i>0$ and all $\alpha\in \Pic0$, that is $E$ satifies I.T. with index $0$.
\end{proof}

The numerical measure of positivity used here is $\Theta$-regularity. Recall from \S2 that
$E$ is $m$-$\Theta$-regular if $E(-(m-1)\Theta)$ is $M$-regular. It is
easy to see that in the case of semihomogeneous bundles this
coincides with strong $m$-$\Theta$-regularity.

\begin{lemma}
A semihomogeneous bundle $E$ is $m$-$\Theta$-regular if and only if $E(-(m-1)\Theta)$ satisfies
I.T. with index $0$.
\end{lemma}
\begin{proof}
The more general fact that an $M$-regular semihomogeneous bundle
satisfies I.T. with index $0$ follows quickly from Lemma
\ref{sh_pullback} above. More precisely the line bundle $M$ in its
statement is forced to be ample since it has a twist with global
sections and positive Euler characteristic.
\end{proof}

Mukai shows in \cite{mukai3} \S6 that the semihomogeneous bundles are Gieseker
semistable (while the simple ones are in fact stable). Moreover, any
semihomogeneous bundle has a Jordan-H\" older filtration in a strong sense:

\begin{proposition}(\cite{mukai3} Proposition 6.18)
Let $E$ be a semihomogeneous bundle on $X$, and let $\delta$ be the equivalence
class of $\frac{{\rm det}(E)}{{\rm rk}(E)}$ in $NS(X)\otimes_\ZZ \QQ$. Then
there exist simple semihomogeneous bundles $F_1, \ldots, F_n$ whose corresponding
class is the same $\delta$, and semihomogeneous bundles $E_1, \ldots, E_n$, satisfying:
\begin{itemize}
\item $E\cong \bigoplus_{i=1}^n E_i$.
\smallskip
\item Each $E_i$ has a filtration whose factors are all isomorphic to $F_i$.
\end{itemize}
\end{proposition}

Since the positivity of $E$ is carried through to the factors of a Jordan-H\"older
filtration as in the Proposition above, standard inductive arguments allow us
to immediately reduce the study below to the case of simple semihomogeneous
bundles, which we do freely in what follows.

\noindent
{\bf A numerical criterion for normal generation.} The main result of this section is that the
normal generation of a semihomogeneous vector bundle is dictated
by an explicit numerical criterion. We assume all throughout
that all the semihomogeneous vector bundles involved satisfy the minimal
positivity condition, namely that they are $0$-$\Theta$-regular (which in particular
is a numerical criterion for global generation, by Theorem \ref{M_reg}).

\begin{theorem}\label{multiplication}
Let $E$ be a rank $r$ semihomogeneous bundle on the polarized
abelian variety $(X,\Theta)$, and assume that $E$ is $0$-$\Theta$-regular.
Then, for any $x\in X$, the multiplication map
$$H^0(E)\otimes H^0(t_x^*E)\longrightarrow H^0(E\otimes t_x^*E)$$
is surjective provided that
$$\frac{1}{r}\cdot c_1(E(-\Theta))+\frac{1}{r^{\prime}}\cdot
\phi_{\Theta}^*c_1(\widehat{E(-\Theta)})>0,$$
where $r^{\prime} = {\rm rk}(\widehat{E(-\Theta)})$
(recall that $\phi_{\Theta}$ is the isogeny $X\rightarrow \hat{X}$ induced
by $\Theta$).
\end{theorem}

\begin{remark}\label{chern}
Although most conveniently written in terms of the Fourier-Mukai transform, the
statement of the theorem is indeed a numerical condition intrinsic
to the vector bundle $E$,
since by \cite{mukai2} Cor. 1.18 one has:
$$c_1(\widehat{E(-\Theta)})=-PD_{2g-2}({\rm ch}_{g-1}(E(-\Theta))),$$
where $PD$ denotes the Poincar\'e duality map
$$PD_{2g-2}:H^{2g-2}(J(X),\ZZ)\rightarrow H^{2}(J(X),\ZZ),$$
and ${\rm ch}_{g-1}$ the $(g-1)$-st
component of the Chern character. Note also that
$${\rm rk}(\widehat{E(-\Theta)})= h^0(E(-\Theta)) = \frac{1}{r^{g-1}}\cdot
\frac{c_1(E(-\Theta))^g}{g!}$$
by \cite{mukai1} Cor. 2.8.
\end{remark}

This implies in a particular case the folowing  explicit statement, one
which will be generalized though at the end of the section.

\begin{corollary}\label{normal_semi}
A $(-1)$-$\Theta$-regular semihomogeneous bundle is normally generated.
\end{corollary}
\begin{proof}
The hypotesis means that $E(-2\Theta)$ satisfies I.T. with index $0$.
Consider an isogeny $f:Y\rightarrow X$ as in Lemma \ref{sh_pushforward},
so that there exists a line bundle $L$ on $Y$ with $f_* L = E(-\Theta)$.
The assumption on $E$ implies that $L$ can be written as the tensor product
of two ample line bundles.
Since $f$ has $0$-dimensional fibers, the Grothendieck-Riemann-Roch theorem
implies that
$$c_1(E(-\Theta))= f_* c_1(L)~{\rm and}~{\rm ch}_{g-1}(E(-\Theta))=
f_* \big(\frac{c_1(L)^{g-1}}{(g-1)!}\big).$$
By Remark \ref{chern} we have then that
$c_1(\widehat{E(-\Theta)})=-PD_{2g-2}\big(f_* \big(\frac{c_1(L)^{g-1}}{(g-1)!}\big)\big)$.
We are then reduced to doing a line bundle computation on $Y$, which follows
by standard methods, as for example in \cite{beauville}.
\end{proof}

\noindent
{\bf Examples.}
There are two basic classes of examples of (strongly) $(-1)$-$\Theta$-regular bundles, and
both turn out to be semihomogeneous. They correspond to the
properties of linear series on abelian varieties and on moduli spaces of
vector bundles on curves, respectively.

\begin{example}({\bf {Projective normality of line bundles.}})
For every ample divisor $\Theta$ on $X$, the line bundle $L
=\OO_X(m\Theta)$ is $(-1)$-$\Theta$-regular for $m\geq 3$. Thus we
recover the classical fact that $\OO_X(m\Theta)$ is projectively
normal for $m\geq 3$.
\end{example}

\smallskip

\begin{example}({\bf {Verlinde bundles.}})\label{Verlinde}
Let $U_C(r,0)$ be the moduli space of rank $r$ and degree $0$ semistable
vector bundles on a
a smooth projective curve $C$ of genus $g\geq 2$. This comes with
a natural determinant map ${\rm det}: U_C(r,0)\rightarrow J(C)$, where
$J(C)$ is the Jacobian of $C$. To a generalized theta
divisor $\Theta_N$ on $U_C(r,0)$ (depending on the choice of a line bundle
$N\in {\rm Pic}^{g-1}(C)$) one associates for any $k\geq 1$
the $(r,k)$-Verlinde bundle
on $J(C)$, defined by $E_{r, k}:= {\rm det}_* \OO(k\Theta_N)$ (cf. \cite{popa}). It is
shown in \emph{loc. cit.} that the numerical properties of $E_{r,k}$ are
essential in understanding the linear series $|k\Theta_N|$ on $U_C(r,0)$.
It is noted there that $E_{r,k}$ are polystable and semihomogeneous.

\noindent
A basic property of these vector bundles is the fact that
$$r_J^*E_{r,k}\cong \oplus \OO_J(kr\Theta_N),$$
where $r_J$ denotes multiplication by $r$ on $J(C)$
(cf. \cite{popa} Lemma 2.3). Noting that the pull-back
$r_J^*\OO_J(\Theta_N)$ is numerically equivalent to
$\OO(r^2\Theta_N)$, we obtain that $E_{r,k}$ is $0$-$\Theta$-regular iff
$k\geq r+1$, and $(-1)$-$\Theta$-regular iff $k\geq 2r+1$. This
implies by the statements above that $E_{r,k}$ is globally generated
for $k\geq r+1$ and normally generated for $k\geq 2r+1$.
These are precisely the results \cite{popa} Proposition 5.2 and
Theorem 5.9(a), the second
obtained there by ad-hoc (though related) methods.
\end{example}

\noindent
{\bf Proof of the numerical criterion.} We will prove a natural generalization
of Theorem \ref{multiplication}, which guarantees the
surjectivity of multiplication maps for two arbitrary semihomogeneous bundles.
This could be seen as an analogue of Butler's theorem \cite{butler} for semistable
vector bundles on curves.

\begin{theorem}\label{mixed_multiplication}
Let $E$ and $F$ be semihomogeneous bundles on $(X,\Theta)$, both satisfying
I.T. with index $0$. Then the multiplication maps
$$H^0(E)\otimes H^0(t_x^*F)\longrightarrow H^0(E\otimes t_x^*F)$$
are surjective for all $x\in X$ if  the following holds:
$$\frac{1}{r_F}\cdot c_1(F(-\Theta)) + \frac{1}{r_E^{\prime}}\cdot
\phi_{\Theta}^* c_1(\widehat{E(-\Theta)})>0.$$
\end{theorem}

We can assume $E$ and $F$ to be simple by the
considerations in \S2, and we will do so in what follows.
We begin with a few technical lemmas.
Let us recall first that, given two sheaves $\E$ and $\G$ on $X$,
their \emph{skew Pontrjagin product} (see \cite{pareschi} \S1) is defined as
$$\E\hat * \G:={p_1}_*((p_1+p_2)^*(\E)\otimes p_2^*(\G)),$$
where $p_1$ and $p_2$ are the projections from $X\times X$ to the two factors.

\begin{lemma}(\cite{pareschi} Theorem 3.1)\label{mult1}
The multiplication map
$$H^0(E)\otimes H^0(t_x^*F)\longrightarrow H^0(E\otimes t_x^*F)$$
is surjective for any $x\in X$ if the skew-Pontrjagin product $E\hat{*}F$ is globally generated,
so in particular if $(E\hat{*} F)$ is $0$-$\Theta$-regular.
\end{lemma}

\begin{lemma}\label{mult2}
For all $i\geq 0$ we have:
$$h^i((E\hat{*} F)\otimes \OO_X(-\Theta)) = h^i((E\hat{*}\OO_X(-\Theta))\otimes F).$$
\end{lemma}
\begin{proof}
This follows from Lemma 3.2 in \cite{pareschi} if we prove the following vanishings:
\begin{enumerate}
\item $h^i(t_x^*E\otimes F) =0,~\forall i>0,~\forall x\in X.$
\item $h^i(t_x^*E\otimes \OO_X(-\Theta))=0,~\forall i>0,~\forall x\in X.$
\end{enumerate}
We treat them separately:
\begin{enumerate}
\item By Lemma \ref{sh_pullback} we know that there exist isogenies $\pi_E:Y_E\rightarrow
X$ and $\pi_F: Y_F\rightarrow X$, and line bundles $M$ on $Y_E$ and $N$ on $Y_F$, such
that $\pi_E^*E\cong \underset{r_E}{\oplus} M$ and $\pi_F^*F\cong \underset{r_F}{\oplus} N$.
Now on the fiber product $Y_E\times_X Y_F$, the pull-back of $t_x^*E\otimes F$ is a direct sum of
line bundles numerically equivalent to $p_1^*M\otimes p_2^*N$ . This line bundle is ample and has sections, and so no
higher cohomology by the Index Theorem. Consequently the same must be true for
$t_x^*E\otimes F$.
\item Since $E$ is semihomogeneous, we have $t_x^*E\cong E\otimes \alpha$ for some
$\alpha\in \Pic0$, and so:
$$h^i(t_x^*E\otimes \OO_X(-\Theta))=h^i(E\otimes\OO_X(-\Theta)\otimes \alpha)=0,$$
since $E(-\Theta)$ satisfies I.T. with index $0$.
\end{enumerate}
\end{proof}

Let us assume from now on for simplicity that the polarization $\Theta$ is symmetric.
This makes the proofs less technical, but the general case is completely similar since
everything depends (via suitable isogenies) only on numerical equivalence classes.

\begin{lemma}\label{mult3}
We have
$$E\hat{*}\OO_X(-\Theta)\cong \phi_{\Theta}^*((-1_X)^*E\otimes \OO_X(-\Theta))^{\widehat{}}
\otimes \OO(-\Theta).$$
\end{lemma}
\begin{proof}
This follows from Mukai's general Lemma 3.10 in \cite{mukai1}.
\end{proof}

Putting together Lemmas \ref{mult1}, \ref{mult2} and \ref{mult3} we get
the following cohomological
criterion for surjectivity of multiplication maps in our given situation.

\begin{proposition}\label{cohomological_multiplication}
Under the hypotheses above, the multiplication maps
$$H^0(E)\otimes H^0(t_x^*F)\longrightarrow H^0(E\otimes t_x^*F)$$
are surjective if we have the following vanishing:
$$h^i(\phi_{\Theta}^*((-1_X)^*E\otimes \OO_X(-\Theta))^{\widehat{}}\otimes F(-\Theta))=0,~\forall i>0.$$
\end{proposition}

We are now in a position to give the proof of Theorem \ref{mixed_multiplication}.
To this end we only need to understand the best numerical assumptions under
which the cohomological requirement in Proposition \ref{cohomological_multiplication}
is satisfied.

\medskip
\begin{proof}(\emph{of Theorem \ref{mixed_multiplication}}.)
We first apply Lemma \ref{sh_pullback} to $G:=\phi_{\Theta}^*\widehat{(-1_X)^*E(-\Theta)}$
and $H:=F(-\Theta)$: there exist isogenies $\pi_G:Y_G\rightarrow X$ and
$\pi_H:Y_H\rightarrow X$, and line bundles $M$ on $Y_G$ and $N$ on $Y_H$, such that
$\pi_G^*G\cong \underset{r_G}{\oplus}M$ and $\pi_H^*H\cong \underset{r_H}{\oplus}N$.
Consider the fiber product $Z:=Y_G\times_X Y_H$, with projections $p_G$ and $p_H$.
Denote by $p:Z\rightarrow X$ the natural composition.

By pulling everything back to $Z$, we see that
$$p^*(G\otimes H)\cong \underset{r_G\cdot r_F}{\bigoplus}(p_1^*M\otimes p_2^*N).$$
This implies that our desired vanishing $H^i(G\otimes H)=0$ (cf. Proposition
\ref{cohomological_multiplication}) holds as long as
$$H^i(p_G^*M\otimes p_H^*N)=0, ~\forall i>0.$$

Now $c_1(p_G^*M)=p_G^*c_1(M)=\frac{1}{r_G}p^*c_1(G)$ and similarly
$c_1(p_H^*N)=p_H^*c_1(N)=\frac{1}{r_H}p^*c_1(G)$.
Finally we get
$$c_1(p_G^*M\otimes p_H^*N)=p^*(\frac{1}{r_G}\cdot c_1(G)+\frac{1}{r_H}\cdot c_1(H)).$$
Thus all we need to have is that the class
$$\frac{1}{r_G}\cdot c_1(G)+\frac{1}{r_H}\cdot c_1(H)$$
be ample, and this is clearly equivalent to the statement of the theorem.
\end{proof}

We conclude by noting that the
reason we only sketched the proof of Corollary \ref{normal_semi} is that in fact under that
particular hypothesis we have
a much more general statement which works for every vector bundle on a polarized abelian variety,
via substantially subtler methods.

\begin{theorem}
Let $E$ and $F$ be $(-1)$-$\Theta$-regular vector bundles on $X$ (i.e.
such that $E(-2\Theta)$ and $F(-2\Theta)$ are $M$-regular).
Then the multiplication map
$$H^0(E)\otimes H^0(F)\rightarrow H^0(E\otimes F)$$
is surjective.
\end{theorem}
\begin{proof}
We use an argument inspired by techniques first introduced by Kempf, and rendered easy
by the results of \cite{pp1}. Let us consider the diagram
$$\xymatrix{
\bigoplus_{\xi\in U} H^0(E(-2\Theta)\otimes P_\xi)\otimes
H^0(2\Theta\otimes P_\xi^\vee)\otimes H^0(F) \ar[r] \ar[d] & H^0(E)\otimes
H^0(F) \ar[d]  \\
\bigoplus_{\xi\in U} H^0(E(-2\Theta)\otimes P_\xi)\otimes H^0(F(2\Theta)\otimes
P_\xi^\vee)  \ar[r] & H^0(E\otimes F) }$$
Under the given hypotheses, the bottom horizontal arrow is onto by the general Theorem \ref{surjectivity}.
On the other hand, the abelian Castelnuovo-Mumford Lemma Theorem \ref{acm} insures that
each one of the components of the vertical map on the left is surjective. Thus the composition is surjective,
which gives the surjectivity of the vertical map on the right.
\end{proof}

\begin{corollary}
Every  $(-1)$-$\Theta$-regular vector bundle is normally generated.
\end{corollary}

\providecommand{\bysame}{\leavevmode\hbox to3em{\hrulefill}\thinspace}

\end{document}